\documentclass{article}
\usepackage{amssymb}
\usepackage{amsfonts}
\usepackage{amsmath}

\begin{document}

\begin{center}
{\Large Applications of some special numbers obtained from a difference
equation of degree three}

\begin{equation*}
\end{equation*}%
Cristina FLAUT and Diana SAVIN 
\begin{equation*}
\end{equation*}
\end{center}

\textbf{Abstract. }{\small In this paper we present applications of special
numbers obtained from a difference equation of degree three. One of these
applications is in the Coding Theory, since some of these numbers can be
used to built cyclic codes with good properties (MDS codes). In another
particular case of these difference equation of degree three, we obtain the
generalized Pell-Fibonacci-Lucas numbers, which were \ extended to the
generalized quaternion algebras. Using properties of these elements, we can
define a set with an interesting algebraic structure, namely an order of a
generalized rational quaternion algebra. }

\bigskip

\medskip


\smallskip

\begin{equation*}
\end{equation*}

\textbf{1. Introduction}%
\begin{equation*}
\end{equation*}

Let $n$ be an arbitrary positive integer and let $a,b,c,x_{0},x_{1},x_{2}~$%
be arbitrary integers. We consider the following difference equation of
degree three%
\begin{equation}
D_{n}=aD_{n-1}+bD_{n-2}+cD_{n-3},D_{0}=x_{0},D_{1}=x_{1},D_{2}=x_{2}. 
\tag{1.1}
\end{equation}%
If we consider $a=b=1,c=0,$ $x_{0}=0,x_{1}=1,x_{2}=1,$ we obtain the
Fibonacci numbers and if we take $a=b=1,c=0,x_{0}=2,x_{1}=1,x_{2}=3,$ we get
the Lucas numbers. If we consider $a=2,b=1,c=0,x_{0}=0,x_{1}=1,x_{2}=2,$ we
obtain the Pell numbers and if we take $a=1,b=0,c=1,x_{0}=0,x_{1}=1,x_{2}=1,$
we find the Fibonacci-Narayana numbers.

Some properties of the above numbers were studied in various paper. In this
paper, we will provide properties and applications of other special numbers
obtained from the equation $\left( 1.1\right) ,$ as from example in the
Coding Theory. Moreover, we extend these numbers to generalized quaternions,
obtaining an interesting algebraic structure. 
\begin{equation*}
\end{equation*}

\begin{equation*}
\end{equation*}

\textbf{2. An application in the Coding Theory}%
\begin{equation*}
\end{equation*}

In [Ba, Pr; 09], [St; 06] and [Ko, Oz, Si; 17] were presented some
applications of the Fibonacci elements in the Coding Theory. In the
following, we will give applications of other special numbers in this domain.

Let $F_{\pi },$ $\pi $ a prime number, be a finite field and let $C\subset
F_{\pi }^{n},$ be a linear code. Let $c=\left(
c_{0},c_{1},c_{2},...,c_{n-1}\right) \in C$ $\ $be a code-word. A linear
code $C$ of length $n$ over a finite field $F_{\pi }~$is a \textit{cyclic
code} if $c=\left( c_{0},c_{1},c_{2},...,c_{n-1}\right) \in C$ implies that $%
c^{\prime }=\left( c_{n-1},c_{0},c_{1},c_{2},...,c_{n-2}\right) \in C.$ From
here, we have that $C$ is invariant at a single right cyclic shift. It is
very useful if each code-word in a cyclic code is represented using
polynomials. Therefore to the code-word $c=\left(
c_{0},c_{1},c_{2},...,c_{n-1}\right) $ we associate the polynomial 
\begin{equation*}
c\left( x\right) =c_{0}+c_{1}x+c_{2}x^{2}+...+c_{n-1}x^{n-1},
\end{equation*}%
called the associated \textit{code polynomial}. It results that the
associated code polynomial for $c^{\prime }$ is 
\begin{equation*}
c^{\prime }\left( x\right)
=c_{n-1}+c_{0}x+c_{1}x^{2}+c_{2}x^{3}+...+c_{n-2}x^{n-1}.
\end{equation*}

We obtain that 
\begin{equation*}
c^{\prime }\left( x\right) =xc\left( x\right) -c_{n-1}\left( x^{n}-1\right) ,
\end{equation*}%
therefore $c^{\prime }\left( x\right) =xc\left( x\right) ~$\textit{mod }$%
\left( x^{n}-1\right) .$ From here, we have that 
\begin{equation*}
c\left( x\right) \in C~mod\left( x^{n}-1\right) \text{ if and only if \ }%
xc\left( x\right) \in C~mod\left( x^{n}-1\right) .
\end{equation*}

A code polynomial $f\left( x\right) $ generates the cyclic code $C$ of
length $n$ if and only if $f\left( x\right) \mid \left( x^{n}-1\right) .$
The polynomial $f$ is called \textit{the generator polynomial} for the code $%
C$. The polynomial $g\left( x\right) =\frac{x^{n}-1}{f\left( x\right) }$ is
called \textit{the check polynomial} of the code \thinspace $C.$

For a code $C,$ the \textit{Hamming distance, }$d,$ between two code-words
of the same length $n~$is the number of positions at which the corresponding
symbols are different. The \textit{minimum} \textit{Hamming distance }of a
code, $d_{H}$ is 
\begin{equation*}
\min \{d_{H}\left( c_{1},c_{2}\right) ~/~c_{1},c_{2}\in C\}.
\end{equation*}
The \textit{Hamming weight, }$w,$ of a code is the number of symbols which
differ from the zero-symbol of the used alphabet. The \textit{minimum} 
\textit{Hamming weight }of a code, $w_{H}$ is 
\begin{equation*}
\min \{w_{H}\left( c_{1},c_{2}\right) ~/~c_{1},c_{2}\in C\}.
\end{equation*}

Let $m\left( x\right) $ be a transmitted message, $r\left( x\right) $ a
received message and $e\left( x\right) $ the error occurred. The \textit{%
syndrome} of the received vector is $s\left( x\right) =r\left( x\right) $ 
\textit{mod}$~f\left( x\right) .$ Supposing that the code $C$ is $t-$errors
corrected code, we compute $s_{1}\left( x\right) =xs\left( x\right) $ 
\textit{mod} $f\left( x\right) $%
\begin{equation*}
s_{1}\left( x\right) =xs\left( x\right) ~\text{\textit{mod} }f\left( x\right)
\end{equation*}%
\begin{equation*}
s_{2}\left( x\right) =x^{2}s\left( x\right) ~\text{\textit{mod} }f\left(
x\right) ....
\end{equation*}

\begin{equation*}
s_{i}\left( x\right) =x^{i}s\left( x\right) ~\text{\textit{mod} }f\left(
x\right) ,
\end{equation*}%
until $w_{H}\left( s_{i}\right) \leq t$. It results the error $e\left(
x\right) =x^{n-i}s_{i}~$\textit{mod} $\left( x^{n}-1\right) ,$ therefore $%
m\left( x\right) =r\left( x\right) -e\left( x\right) .$For \ other details
regarding cyclic codes, the reader is referred, for example, to [Li, Xi;
04].\medskip

\textbf{Definition 2.1. }Let \thinspace $C$ be a linear code, with
parameters $[n,k,d_{H}]$, with $n$ the length of the codewords, $k$ the
dimension of the code$.$ If we have $k+d=n+1,$ the code $C$ is called 
\textit{maximum distance separable} code, shortly MDS.\medskip

Let $\pi $ be a positive prime integer. We consider the degree three
sequence given in $\left( 1.1\right) .$ We denote with $l_{d}\left( \pi
\right) $ and $\beta _{d}\left( \pi \right) $ the period, respectively the
number of zeros in a single period for this sequence defined on $\mathbb{Z}%
_{\pi },$where $\mathbb{Z}_{\pi }$ denotes the finite field of integers
modulo $\pi ,$ having $\pi $ elements$.\medskip $

We consider $\delta \left( x\right) \in \mathbb{Z}_{\pi }[x],\delta \left(
x\right) =\overset{l_{d}\left( \pi \right) }{\underset{i=0}{\sum }}%
D_{i}x^{i},$ the $\mathbf{D}-$\textit{polynomial} associated to the sequence 
$\mathbf{D,}$ given by $\left( 1.1\right) ,$ over \ $\mathbb{Z}_{\pi },$
where $\mathbf{D}=\{D_{0},D_{1},...D_{n},...\}.\medskip $

\textbf{Theorem 2.2. }\textit{With the above notations, let} $\delta \left(
x\right) \in \mathbb{Z}_{\pi }[x]$ \textit{be the }$\mathbf{D}-$\textit{\
polynomial. The following relation is true}%
\begin{equation*}
\delta \left( x\right) \left( cx^{3}+bx^{2}+ax-1\right) =
\end{equation*}%
\begin{equation*}
=cD_{l-1}x^{l+2}+D_{l+1}x^{l+1}+\left( aD_{1}+bD_{0}-D_{2}\right)
x^{2}+\left( aD_{0}-D_{1}\right) x-D_{0}.
\end{equation*}%
$.$

\textbf{Proof.} We denote $l_{D}\left( \pi \right) =l.$ Let $%
D_{n}=aD_{n-1}+bD_{n-2}+cD_{n-3},d_{0}=x_{0},d_{1}=x_{1},d_{2}=x_{2},$be a
difference equation of degree three. First of all, we compute\newline
$\delta \left( x\right) \left( cx^{3}+bx^{2}+ax-1\right) =\left(
cx^{3}+bx^{2}+ax-1\right) \overset{l-1}{\underset{i=0}{\sum }}D_{i}x^{i}=$%
\newline
$=cD_{l-1}x^{l+2}+bD_{l-1}x^{l+1}+aD_{l-1}x^{l}-D_{l-1}x^{l-1}+$\newline
$+cD_{l-2}x^{l+1}+bD_{l-2}x^{l}+aD_{l-2}x^{l-1}-D_{l-2}x^{l-2}+$\newline
$+cD_{l-3}x^{l}+bD_{l-3}x^{l-1}+aD_{l-3}x^{l-2}-D_{l-3}x^{l-3}+$\newline
$+cD_{l-4}x^{l-1}+bD_{l-4}x^{l-2}+aD_{l-4}x^{l-3}-D_{l-4}x^{l-4}+...+$%
\newline
$+cD_{4}x^{7}+bD_{4}x^{6}+aD_{4}x^{5}-D_{4}x^{4}+$\newline
$+cD_{3}x^{6}+bD_{3}x^{5}+aD_{3}x^{4}-D_{3}x^{3}+$\newline
$+cD_{2}x^{5}+bD_{2}x^{4}+aD_{2}x^{3}-D_{2}x^{2}+$\newline
$+cD_{1}x^{4}+bD_{1}x^{3}+aD_{1}x^{2}-D_{1}x+$\newline
$+cD_{0}x^{3}+bD_{0}x^{2}+aD_{0}x-D_{0}.$\newline
Using relation $\left( 1.1\right) $, we obtain\newline
$\delta \left( x\right) \left( cx^{3}+bx^{2}+ax-1\right) =$\newline
$=cD_{l-1}x^{l+2}+bD_{l-1}x^{l+1}+aD_{l-1}x^{l}+$\newline
$+cD_{l-2}x^{l+1}+bD_{l-2}x^{l}+cD_{l-3}x^{l}+$\newline
$+aD_{1}x^{2}-D_{1}x-D_{2}x^{2}+bD_{0}x^{2}+aD_{0}x-D_{0}=$\newline
$=cD_{l-1}x^{l+2}+D_{l+1}x^{l+1}+\left( aD_{1}+bD_{0}-D_{2}\right)
x^{2}+\left( aD_{0}-D_{1}\right) x-D_{0}.\Box \smallskip $

\textbf{Remark 2.3.} For $D_{0}=0,D\,_{1}=1,c=0,b=1,$ since $D_{2}=aD_{1}$ $%
D_{l+1}=D_{1}~$\textit{mod} $\pi ,$ we have $\delta \left( x\right) \left(
cx^{3}+bx^{2}+ax-1\right) =\delta \left( x\right) \left( x^{2}+ax-1\right) =$%
\newline
$=D_{l+1}x^{l+1}+\left( aD_{1}-D_{2}\right) x^{2}-x=x^{l+1}-x.~$Therefore,
in this situation, we get $\delta \left( x\right) \left( x^{2}+ax-1\right)
=x^{l+1}-x.\medskip $

We consider the following difference equation 
\begin{equation}
D_{n}=aD_{n-1}+bD_{n-2},  \tag{2.1}
\end{equation}%
with $D_{0}=0,D\,_{1}=1.\medskip $

\textbf{Theorem 2.4. }\textit{With the above notations, let} $\delta \left(
x\right) \in \mathbb{Z}_{\pi }[x]$ \textit{be the }$\mathbf{D}-$\textit{%
polynomial associated to the sequence }$\left( 2.1\right) ,~\delta \left(
x\right) =x\delta ^{\prime }\left( x\right) $\textit{. The cyclic code} \ $%
C=<\delta ^{\prime }\left( x\right) >,$ \textit{generated by} $\delta
^{\prime }\left( x\right) ,$ \textit{has dimension} $2$ \textit{and minimum
Hamming distance} $d=l_{d}\left( \pi \right) -\beta _{d}\left( \pi \right)
.\medskip $

\textbf{Proof.} We denote $l_{d}\left( \pi \right) =l.$ Using Theorem 2.2
and Remark 2.3, for $D_{0}=0,D\,_{1}=1,c=0,b=1,$we obtain that\newline
$\delta \left( x\right) \left( x^{2}+2x-1\right) =x^{l+1}-x=x\left(
x^{l}-1\right) .$ It is clear that $\delta ^{\prime }\left( x\right) \mid
\left( x^{l}-1\right) ,$ $x\mid \delta \left( x\right) $ and $x\nmid \left(
x^{l}-1\right) .$ From here, we obtain that gcd$\left( \delta \left(
x\right) ,x^{l}-1\right) =\frac{x^{l}-1}{x^{2}+ax-1}=\frac{\delta \left(
x\right) }{x}=\delta ^{\prime }\left( x\right) .$ We consider the cyclic
code of length $l,~$generated by $\delta ^{\prime }\left( x\right)
,~C=<\delta ^{\prime }\left( x\right) >.$ This code has dimension $l-~$deg$%
\left( \delta ^{\prime }\left( x\right) \right) =2$, therefore has $\pi ^{2}$
elements. From above, $\delta ^{\prime }\left( x\right) =\frac{\delta \left(
x\right) }{x},$ therefore the polynomials $\delta ^{\prime }\left( x\right) $
and $\delta \left( x\right) $ have the same weight. Since the number of the
zero elements in a single period for a $D-$sequence is $\beta _{d}\left( \pi
\right) ,$ we have that the minimum Hamming weight is $w\left( \delta
^{\prime }\left( x\right) \right) =l-\beta _{d}\left( \pi \right) .$ To
prove that the minimum Hamming distance is $l-\beta _{d}\left( \pi \right) ,$
we consider the polynomials $\delta ^{\prime }\left( x\right) $ and $\delta
\left( x\right) =x\delta ^{\prime }\left( x\right) \in C.$ $\ $The
polynomial $\delta ^{\prime }\left( x\right) $ has $\beta _{d}\left( \pi
\right) ~$zero elements on the positions translated with one to the right,
compared with the position of the zeros in the given period \ and the
polynomial $x\delta ^{\prime }\left( x\right) $ has $\beta _{d}\left( \pi
\right) ~$ zero elements on the positions given by the period. If we
consider another polynomial $c\left( x\right) \in $ $C$ with the number of
zero elements more than $\beta _{d}\left( \pi \right) ,$ supposing for
example $\beta _{d}\left( \pi \right) +1,$ we have that $\delta ^{\prime
}\left( x\right) ,x\delta ^{\prime }\left( x\right) ,c\left( x\right) $ are
linearly independent elements, therefore the minimum Hamming distance is $%
l-\beta _{d}\left( \pi \right) .\Box \smallskip $

The above Theorem generalized Theorems 2.2, 2.4, 2.6 from [Ko, Oz, Si; 17]
to $\mathbf{D}-$polynomials. Using ideas from the above mentioned paper,
where were presented cyclic codes obtained from Fibonacci polynomials, in
the following, for particular cases, we will present some properties of
cyclic codes obtained from $D-$polynomials.\qquad

Let $\pi $ be a positive prime integer. We denote with $l_{fib}\left( \pi
\right) $ and $\beta _{fib}\left( \pi \right) $ the period, respectively the
number of zeros in a single period of the Fibonacci sequence considered in $%
\mathbb{Z}_{\pi }$ and with $l_{pell}\left( \pi \right) $ and $\beta
_{fib}\left( \pi \right) $ \ the period, respectively the number of zeros in
a single period of of the Pell sequence considered in $\mathbb{Z}_{\pi
}.\medskip \medskip $\newline

\textbf{Example 2.5. \ }Let $D_{0}=0,D\,_{1}=1,c=0,b=1,$

i) For $a=1,$ we have the Fibonacci numbers. For $\pi =3,$ we have $%
l_{fib}\left( \pi \right) =8,\beta _{fib}\left( \pi \right) =2$, for $\pi
=5, $ we get $l_{fib}\left( \pi \right) =20,\beta _{fib}\left( \pi \right)
=4,$ for $\pi =7,$ we have $l_{fib}\left( \pi \right) =16,\beta _{fib}\left(
\pi \right) =2,$ for $\pi =11,$ we obtain $l_{fib}\left( \pi \right)
=10,\beta _{fib}\left( \pi \right) =1,$ for $\pi =13,$ we have $%
l_{fib}\left( \pi \right) =28,\beta _{fib}\left( \pi \right) =4.$ (Also you
can see [Ko, Oz, Si; 17])

ii) For $a=2,$ we get the Pell numbers. For $\pi =3,$ we have $%
l_{pell}\left( \pi \right) =8,\beta _{pell}\left( \pi \right) =2.$ If $\pi
=5,$ we obtain $l_{pell}\left( \pi \right) =12,\beta _{pell}\left( \pi
\right) =4$ and if $\pi =7,$ we have $l_{pell}\left( \pi \right) =6,\beta
_{pell}\left( \pi \right) =1.$ For $\pi =11,$ we have $l_{pell}\left( \pi
\right) =24,\beta _{pell}\left( \pi \right) =1$ and for $\pi =13,$ we get $%
l_{pell}\left( \pi \right) =28,\beta _{pell}\left( \pi \right) =4.$

iii) Case \thinspace $a=3.$ For example, if $\pi =11,$ we have $l_{3}\left(
\pi \right) =8,\beta _{3}\left( \pi \right) =2.$

iv) Case \thinspace $a=4.$ For example, if $\pi =11,$ we get $l_{4}\left(
\pi \right) =10,\beta _{4}\left( \pi \right) =1.$\ 

v) Case \thinspace $a=5.$ For example, if $\pi =7,$ we obtain $l_{5}\left(
\pi \right) =6,\beta _{5}\left( \pi \right) =1.$

vi) Case \thinspace $a=6.$\textbf{\ }For example, if $\pi =7,$ we have $%
l_{7}\left( \pi \right) =4,\beta _{7}\left( \pi \right) =1~$and if $\pi =13,$
we get $l_{13}\left( \pi \right) =6,\beta _{13}\left( \pi \right)
=1.\medskip $

\textbf{Example 2.6.} i) We consider the Pell numbers for $\pi =7.$ We have $%
l_{pell}\left( \pi \right) =6,\beta _{pell}\left( \pi \right) =1.~$The
corresponding Pell sequence modulo $7$ is\newline
$0,1,2,5,5,1,$ and we obtain the Pell-polynomial, $p\left( x\right)
=x+2x^{2}+5x^{3}+5x^{4}+x^{5}.$ The code $C,$ generated by the polynomial $%
p^{\prime }\left( x\right) ,$ has minimum Hamming distance $d=5,$ and it is
on the type $[6,2,5].$Since $k+d=7=n+1,$ the code $C$ is an MDS code.

ii) We consider $a=5$ and $\pi =7.$ We have $l_{5}\left( \pi \right)
=6,\beta _{5}\left( \pi \right) =1.$The corresponding $\mathbf{D}-$sequence
modulo $7$ is $0,1,5,5,2,1$ and we obtain the $\mathbf{D}$-polynomial $%
\delta \left( x\right) =x+5x^{2}+5x^{3}+2x^{4}+x^{5}.$ The code $C,$
generated by the polynomial $\delta ^{\prime }\left( x\right) ,$ has Hamming
distance $d=5~$and it is on the type $[6,2,5],$ being an MDS code.

iii) We consider $a=6$ and $\pi =13.$ We have $l_{6}\left( \pi \right)
=6,\beta _{6}\left( \pi \right) =1.$The corresponding $\mathbf{D}-$sequence
modulo $13$ is $0,1,6,11,7,1$ and we obtain the $\mathbf{D}$-polynomial $%
\delta \left( x\right) =x+6x^{2}+11x^{3}+7x^{4}+x^{5}.$ The code $C,$
generated by the polynomial $\delta ^{\prime }\left( x\right) ,$ has Hamming
distance $d=5,$ and it is on the type $[6,2,5],$ being an MDS code.\medskip

\textbf{Example 2.7.} i) We consider $a=3$ and $\pi =11.$ We have $%
l_{3}\left( \pi \right) =8,\beta _{3}\left( \pi \right) =2.$The
corresponding $\mathbf{D}-$sequence modulo $11$ is $0,1,3,10,0,10,8,1$ and
we obtain the $\mathbf{D}$-polynomial $\delta \left( x\right)
=x+3x^{2}+10x^{3}+10x^{5}+8x^{6}+x^{7}.$ The code $C,$ generated by the
polynomial $\delta ^{\prime }\left( x\right) ,$ has Hamming distance $d=6,$
and it is on the type $[8,2,6].~$Since $k+d=7=n+1,$ this code is not an MDS
code.

ii) For $a=8$ and $\pi =13,$ we have $l_{8}\left( \pi \right) =12,\beta
_{8}\left( \pi \right) =4.$The corresponding $\mathbf{D}-$sequence modulo $%
13 $ is $0,1,8,0,8,12,0,12,5,0,5,1$ and we obtain the $\mathbf{D}$%
-polynomial $\delta \left( x\right)
=x+8x^{2}+8x^{4}+12x^{5}+12x^{7}+5x^{8}+5x^{10}+x^{11}.$ The code $C,$
generated by the polynomial $\delta ^{\prime }\left( x\right) ,$ has Hamming
distance $d=8,$ is on the type $[8,2,8]$ and it is not an MDS code.\medskip

\textbf{Example 2.8.} We consider for $\pi =7$, the Pell-polynomial, $%
p\left( x\right) =x^{5}+5x^{4}+5x^{3}+$ $2x^{2}+x~$and the code $C$
generated by the polynomial $p^{\prime }\left( x\right)
=x^{4}+5x^{3}+5x^{2}+2x+1.~$This code is $2-$errors correcting code. We
suppose that the code-word sent was $c=\left( 0,3,6,1,1,3\right) $ and we
received the code-word $r=\left( 0,5,6,1,1,6\right) .$ The corresponding
polynomial for the code-word $c$ is $3x+6x^{2}+x^{3}+x^{4}+3x^{5}.$ We
compute the syndrome polynomial. We get\newline
$6x^{5}+x^{4}+x^{3}+6x^{2}+5x=$\newline
$=\left( x^{4}+5x^{3}+5x^{2}+2x+1\right) \left( 6x+6\right)
+4x^{3}+6x^{2}+x+1,$ therefore the syndrome is $s\left( x\right)
=(4x^{3}+6x^{2}+x+1)$ \textit{mod} $\left( x^{4}+5x^{3}+5x^{2}+2x+\right) .$
We have\newline
$s_{1}\left( x\right) =xs\left( x\right) $~\textit{mod} $\left(
x^{4}+5x^{3}+5x^{2}+2x+1\right) =$\newline
$=(4x^{4}+6x^{3}+x^{2}+x)$ \textit{mod} $\left(
x^{4}+5x^{3}+5x^{2}+2x+1\right) =$\newline
$=4\left( x^{4}+5x^{3}+5x^{2}+2x+1\right) +2x^{2}+3~$\textit{mod} $\left(
x^{4}+5x^{3}+5x^{2}+2x+1\right) .$ \newline
Therefore $s_{1}\left( x\right) =2x^{2}+3$ \textit{mod} $p^{\prime }.\,\ $It
results that $i=1$ and the error polynomial is $e\left( x\right)
=x^{6-1}s_{1}\left( x\right) $ \textit{mod }$\left( x^{6}-1\right) .$ Since $%
2x^{7}+3x^{5}=2x\left( x^{6}-1\right) +3x^{5}+2x,$ we obtain that $%
(2x^{7}+3x^{5})~$\textit{mod }$\left( x^{6}-1\right) =3x^{5}+2x.$ Therefore,
the submitted word was $c=r-e=$ $\left( 0,5,6,1,1,6\right) -\left(
0,2,0,0,0,3\right) =\left( 0,3,6,1,1,3\right) ,$that means we recover the
submitted word.%
\begin{equation*}
\end{equation*}

\textbf{3.} \textbf{Some applications of} \textbf{generalized
Pell-Fibonacci-Lucas elements } 
\begin{equation*}
\end{equation*}

$\allowbreak $In relation $\left( 1.1\right) ,$ if we consider $%
a=2,b=1,c=0,x_{0}=0,x_{1}=1,x_{2}=2,$ we obtain the Pell numbers and if we
take $a=2,b=1,c=0,x_{0}=2,x_{1}=2,x_{2}=6,$ we obtain the Pell-Lucas
numbers. Let $\left( P_{n}\right) _{n\geq 0}$ be the Pell sequence 
\begin{equation*}
P_{n}=2P_{n-1}+P_{n-2},\;n\geq 2,P_{0}=0;P_{1}=1,
\end{equation*}%
and $\left( Q_{n}\right) _{n\geq 0}$ be the Pell-Lucas sequence 
\begin{equation*}
Q_{n}=2Q_{n-1}+Q_{n-2},\;n\geq 2,Q_{0}=2;Q_{1}=2.
\end{equation*}%
We consider the numbers $\alpha =1+\sqrt{2}$ and $\beta =1-\sqrt{2}.$ The
following formulae are well known:\newline
\smallskip \newline
\textbf{Binet's formula for Pell sequence.} 
\begin{equation*}
P_{n}=\frac{\alpha ^{n}-\beta ^{n}}{\alpha -\beta }=\frac{\alpha ^{n}-\beta
^{n}}{2\sqrt{2}},\ \ \left( \forall \right) n\in \mathbb{N}.
\end{equation*}%
\textbf{Binet's formula for Pell-Lucas sequence.} 
\begin{equation*}
Q_{n}=\alpha ^{n}+\beta ^{n},\ \ \left( \forall \right) n\in \mathbb{N}.
\end{equation*}

\medskip

Let $A$ be the generating function for the sequence $\left( D_{n}\right)
_{n\geq 0}$, \ given by the relation $\left( 1,1\right) $ $A\left( z\right)
=\sum\limits_{n\geq 0}^{{}}D_{n}z^{n}.$ In the following, we determine this
function.\newline

\textbf{Proposition 3.1.} \textit{We have:} 
\begin{equation*}
A\left( z\right) =\frac{D_{0}+\left( D_{1}-aD_{0}\right) z+\left(
D_{2}-aD_{1}-bD_{0}\right) z^{2}}{1-az-bz^{2}-cz^{3}}.
\end{equation*}%
\textbf{Proof.} 
\begin{equation}
A\left( z\right) =D_{0}+D_{1}z+D_{2}z^{2}+D_{3}z^{3}+...+D_{n}z^{n}+... 
\tag{3.1}
\end{equation}%
\begin{equation}
azA\left( z\right) =aD_{0}z+aD_{1}z^{2}+aD_{2}z^{3}+...+aD_{n-1}z^{n}+... 
\tag{3.2}
\end{equation}%
\begin{equation}
bz^{2}A\left( z\right)
=bD_{0}z^{2}+bD_{1}z^{3}+bD_{2}z^{4}+...+bD_{n-2}z^{n}+...  \tag{3.3}
\end{equation}%
\begin{equation}
cz^{3}A\left( z\right)
=cD_{0}z^{3}+cD_{1}z^{4}+cD_{2}z^{5}+...+cD_{n-3}z^{n}+...  \tag{3.4}
\end{equation}%
Adding the equalities (3.2), (3.3) and (3.4) member by member, we obtain: 
\begin{equation*}
A\left( z\right) \left( 1-az-bz^{2}-cz^{3}\right) =D_{0}+\left(
D_{1}-aD_{0}\right) z+\left( D_{2}-aD_{1}-bD_{0}\right) z^{2}.
\end{equation*}%
Therefore, we get%
\begin{equation*}
A\left( z\right) =\frac{D_{0}+\left( D_{1}-aD_{0}\right)
z+(D_{2}-aD_{1}-bD_{0})z^{2}}{1-az-bz^{2}-cz^{3}}.
\end{equation*}%
$\Box \smallskip $\medskip \newline
For $a=c=1$ and $b=3$ we obtain the sequence $\left( D_{n}\right) _{n\geq
0}, $ $D_{0}=0,$ $D_{1}=D_{2}=1,$ $D_{n+3}=D_{n+2}+3D_{n+1}+D_{n},\;n\geq 0.$
This equality is equivalent with 
\begin{equation*}
D_{n+3}+D_{n+2}=2\left( D_{n+2}+D_{n+1}\right) +\left( D_{n+1}+D_{n}\right) .
\end{equation*}%
If we take the sequence $\left( b_{n}\right) _{n\geq 0},$ $%
b_{n+1}=D_{n+1}+D_{n},\;n\geq 0,$ the last equality becomes 
\begin{equation*}
b_{n+3}=2b_{n+2}+b_{n+1},\;n\geq 0,
\end{equation*}%
where $b_{1}=1$ and $b_{2}=2.$ Moreover, if we consider $b_{0}=0,$ it
results that the sequence $\left( b_{n}\right) _{n\geq 0}$ is in fact the
sequence of Pell numbers $\left( P_{n}\right) _{n\geq 0}.\medskip $

\textbf{Proposition 3.2.} \textit{Let} $\left( P_{n}\right) _{n\geq 0}$ 
\textit{be the sequence of Pell numbers and} $\left( Q_{n}\right) _{n\geq 0}$
\textit{be the sequence of Pell-Lucas numbers. Let} $A$ \textit{be the
following matrix}\newline
$A=\left( 
\begin{array}{lll}
\frac{Q_{1}}{2} & 0 & \sqrt{2}P_{1} \\ 
0 & \frac{Q_{1}}{2}+\sqrt{2}P_{1} & 0 \\ 
\sqrt{2}P_{1} & 0 & \frac{Q_{1}}{2}%
\end{array}%
\right) .$ \textit{Then, we have}: 
\begin{equation*}
A^{n}=\left( 
\begin{array}{lll}
\frac{Q_{n}}{2} & 0 & \sqrt{2}P_{n} \\ 
0 & \frac{Q_{n}}{2}+\sqrt{2}P_{n} & 0 \\ 
\sqrt{2}P_{n} & 0 & \frac{Q_{n}}{2}%
\end{array}%
\right) .
\end{equation*}%
\textbf{Proof.} We prove by induction after $n$$\in $$\mathbb{N}^{\ast }$
the following statement\newline
\begin{equation*}
P\left( n\right) :A^{n}=\left( 
\begin{array}{lll}
\frac{Q_{n}}{2} & 0 & \sqrt{2}P_{n} \\ 
0 & \frac{Q_{n}}{2}+\sqrt{2}P_{n} & 0 \\ 
\sqrt{2}P_{n} & 0 & \frac{Q_{n}}{2}%
\end{array}%
\right) .
\end{equation*}%
We remark that $P\left( 1\right) $ is true.\newline
We suppose that $P\left( n\right) $ is true and we prove that $P\left(
n+1\right) $ is true. 
\begin{equation*}
P\left( n+1\right) :A^{n+1}=\left( 
\begin{array}{lll}
\frac{Q_{n+1}}{2} & 0 & \sqrt{2}P_{n+1} \\ 
0 & \frac{Q_{n+1}}{2}+\sqrt{2}P_{n+1} & 0 \\ 
\sqrt{2}P_{n+1} & 0 & \frac{Q_{n+1}}{2}%
\end{array}%
\right) .
\end{equation*}%
\begin{equation*}
A^{n+1}=A^{n}A=\left( 
\begin{array}{lll}
\frac{Q_{n}}{2} & 0 & \sqrt{2}P_{n} \\ 
0 & \frac{Q_{n}}{2}+\sqrt{2}P_{n} & 0 \\ 
\sqrt{2}P_{n} & 0 & \frac{Q_{n}}{2}%
\end{array}%
\right) \cdot
\end{equation*}%
\begin{equation*}
\cdot \left( 
\begin{array}{lll}
\frac{Q_{1}}{2} & 0 & \sqrt{2}P_{1} \\ 
0 & \frac{Q_{1}}{2}+\sqrt{2}P_{1} & 0 \\ 
\sqrt{2}P_{1} & 0 & \frac{Q_{1}}{2}%
\end{array}%
\right).
\end{equation*}%
Using Binet's formulas for Pell-Lucas sequence and Pell-Lucas sequence, we
obtain: 
\begin{equation*}
A^{n+1}=\left( 
\begin{array}{lll}
\frac{\alpha ^{n}+\beta ^{n}}{2} & 0 & \frac{\alpha ^{n}-\beta ^{n}}{2} \\ 
0 & \alpha ^{n} & 0 \\ 
\frac{\alpha ^{n}-\beta ^{n}}{2} & 0 & \frac{\alpha ^{n}+\beta ^{n}}{2}%
\end{array}%
\right) \cdot \left( 
\begin{array}{lll}
\frac{\alpha +\beta }{2} & 0 & \frac{\alpha -\beta }{2} \\ 
0 & \alpha & 0 \\ 
\frac{\alpha -\beta }{2} & 0 & \frac{\alpha +\beta }{2}%
\end{array}%
\right) =
\end{equation*}%
\begin{equation*}
\left( 
\begin{array}{lll}
\frac{\alpha ^{n+1}+\beta ^{n+1}}{2} & 0 & \frac{\alpha ^{n+1}-\beta ^{n+1}}{%
2} \\ 
0 & \alpha ^{n+1} & 0 \\ 
\frac{\alpha ^{n+1}-\beta ^{n+1}}{2} & 0 & \frac{\alpha ^{n+1}+\beta ^{n+1}}{%
2}%
\end{array}%
\right) =\left( 
\begin{array}{lll}
\frac{Q_{n+1}}{2} & 0 & \sqrt{2}P_{n+1} \\ 
0 & \frac{Q_{n+1}}{2}+\sqrt{2}P_{n+1} & 0 \\ 
\sqrt{2}P_{n+1} & 0 & \frac{Q_{n+1}}{2}%
\end{array}%
\right) .
\end{equation*}%
Therefore, $P\left( n+1\right) $ is true.$\Box \smallskip $

In the paper [Fl, Sa; 15], we introduced the generalized Fibonacci-Lucas
numbers and the generalized Fibonacci-Lucas quaternions and we obtained many
properties of these elements ( also you can see [Sa; 17]). In a similar way,
we introduce here the generalized Pell-Fibonacci-Lucas numbers and the
generalized Pell-Fibonacci-Lucas quaternions.\newline
First, we will give some identities involving Pell numbers and Pell-Lucas
numbers.\newline
\medskip \newline
\textbf{Proposition 3.3. } \textit{Let} $\left( P_{n}\right) _{n\geq 0}$ 
\textit{be the sequence of Pell numbers and} $\left( Q_{n}\right) _{n\geq 0}$
\textit{be the sequence of Pell-Lucas numbers. Then, we have}:\newline
i) 
\begin{equation*}
Q_{n}Q_{n+l}=Q_{2n+l}+\left( -1\right) ^{n}Q_{l},\ n,l\in \mathbb{N};
\end{equation*}%
ii) 
\begin{equation*}
P_{n}Q_{n+l}=P_{2n+l}+\left( -1\right) ^{n+1}P_{l},\ n,l\in \mathbb{N};
\end{equation*}%
iii) 
\begin{equation*}
P_{n+l}Q_{n}=P_{2n+l}+\left( -1\right) ^{n}P_{l},\ n\in \mathbb{N};
\end{equation*}%
iv) 
\begin{equation*}
P_{n}P_{n+l}=\frac{1}{8}\left( Q_{2n+l}+\left( -1\right) ^{n+1}Q_{l}\right)
,~n,l\in \mathbb{N}.
\end{equation*}%
\textbf{Proof.} We will use Binet's formulas for Pell numbers and Pell-Lucas
numbers.\newline
i) 
\begin{equation*}
Q_{n}Q_{n+l}=\left( \alpha ^{n}+\beta ^{n}\right) \left( \alpha ^{n+l}+\beta
^{n+l}\right) =\alpha ^{2n+l}+\beta ^{2n+l}+\alpha ^{n}\beta ^{n+l}+\alpha
^{n+l}\beta ^{n}=
\end{equation*}%
\begin{equation*}
=Q_{2n+l}+\alpha ^{n}\beta ^{n}Q_{l}=Q_{2n+l}+\left( -1\right) ^{n}Q_{l}.
\end{equation*}%
ii) 
\begin{equation*}
P_{n}Q_{n+l}=\frac{\alpha ^{n}-\beta ^{n}}{\alpha -\beta }\left( \alpha
^{n+l}+\beta ^{n+l}\right) =
\end{equation*}%
\begin{equation*}
=\frac{\alpha ^{2n+l}-\beta ^{2n+l}}{\alpha -\beta }-\frac{\alpha ^{n}\beta
^{n}\left( \alpha ^{l}-\beta ^{l}\right) }{\alpha -\beta }=P_{2n+l}+\left(
-1\right) ^{n+1}P_{l}.
\end{equation*}%
iii) 
\begin{equation*}
P_{n+l}Q_{n}=\frac{\alpha ^{n+l}-\beta ^{n+l}}{\alpha -\beta }\left( \alpha
^{n}+\beta ^{n}\right) =
\end{equation*}%
\begin{equation*}
=\frac{\alpha ^{2n+l}-\beta ^{2n+l}}{\alpha -\beta }+\frac{\alpha ^{n}\beta
^{n}\left( \alpha ^{l}-\beta ^{l}\right) }{\alpha -\beta }=P_{2n+l}+\left(
-1\right) ^{n}P_{l}.
\end{equation*}%
iv) 
\begin{equation*}
P_{n}P_{n+l}=\frac{\alpha ^{n}-\beta ^{n}}{\alpha -\beta }\cdot \frac{\alpha
^{n+l}-\beta ^{n+l}}{\alpha -\beta }=
\end{equation*}%
\begin{equation*}
=\frac{\alpha ^{2n+l}+\beta ^{2n+l}-\alpha ^{n}\beta ^{n}\left( \alpha
^{l}+\beta ^{l}\right) }{8}=\frac{1}{8}\left( Q_{2n+l}+\left( -1\right)
^{n+1}Q_{l}\right) .
\end{equation*}%
$\Box \smallskip $\medskip \newline

Let $n$ be an arbitrary positive integer and $p,q$ be two arbitrary
integers. We introduce the following numbers $\left( r_{n}\right) _{n\geq 1},
$%
\begin{equation*}
r_{n+1}=pP_{n}+qQ_{n+1},\;n\geq 0.
\end{equation*}%
We call the numbers $r_{n}$ $(n\geq 1)$ the \textit{generalized
Pell-Fibonacci-Lucas numbers}.\newline
To avoid confusions, in the followings we will use instead of $r_{n}$ the
notation $r_{n}^{p,q}.$ \smallskip

In the following, we consider $\alpha,\beta\in \mathbb{Q}\setminus \{0\}$.
Let $\mathbb{H}\left( \alpha,\beta\right) $ be the generalized quaternion%
{\small \ }algebra with basis $\{1,e_{1},e_{2},e_{3}\}.$ This algebra is a
rational algebra with the multiplication given in the following table%
\begin{equation*}
\begin{tabular}{c||c|c|c|c|}
$\cdot $ & $1$ & $e_{1}$ & $e_{2}$ & $e_{3}$ \\ \hline\hline
$1$ & $1$ & $e_{1}$ & $e_{2}$ & $e_{3}$ \\ \hline
$e_{1}$ & $e_{1}$ & $\alpha$ & $e_{3}$ & $\alpha e_{2}$ \\ \hline
$e_{2}$ & $e_{2}$ & $-e_{3}$ & $\beta$ & $-\beta e_{1}$ \\ \hline
$e_{3}$ & $e_{3}$ & $-\alpha e_{2}$ & $\beta e_{1}$ & $\gamma _{1}\beta$ \\ 
\hline
\end{tabular}%
.
\end{equation*}%
We consider $\alpha ,\beta $$\in $$\mathbb{Q}^{\ast }$ and $\mathbb{H}_{%
\mathbb{Q}}\left( \alpha ,\beta \right) $ the generalized quaternion algebra
over the rational field$.$ We define the $n$-th \textit{generalized Pell-
Fibonacci-Lucas quaternion }to be the element 
\begin{equation*}
R_{n}^{p,q}=r_{n}^{p,q}\cdot 1+r_{n+1}^{p,q}i+r_{n+2}^{p,q}j+r_{n+3}^{p,q}k.
\end{equation*}

\textbf{Remark 3.4. } \textit{Let} $n$ \textit{be an arbitrary positive
integer and} $p,q$ \textit{be two arbitrary integers. Let }$\left(
r_{n}^{p,q}\right) _{n\geq 1}$ \textit{be the generalized
Pell-Fibonacci-Lucas numbers and} $\left( R_{n}^{p,q}\right) _{n\geq 1}$ 
\textit{be the generalized Pell-Fibonacci-Lucas quaternion elements. Then we
have the following relation} 
\begin{equation*}
R_{n}^{p,q}=0\ \text{\textit{if\ and\ only\ if}}\ p=q=0.
\end{equation*}%
\textbf{Proof.} " $\Leftarrow $" It is trivial.\newline
" $\Rightarrow $" Using the fact that $\left\{ 1,i,j,k\right\} $ is a basis
in $\mathbb{H}_{\mathbb{Q}}\left( \alpha ,\beta \right) ,$ we obtain that $%
r_{n}^{p,q}=r_{n+1}^{p,q}=r_{n+2}^{p,q}=r_{n+3}^{p,q}=0.$ This implies that $%
r_{n-1}^{p,q}=0$ ..., $r_{2}^{p,q}=0$, $r_{1}^{p,q}=0.$ Since $%
r_{1}^{p,q}=pP_{0}+qQ_{1}=2q,$ it results $q=0.$ From $r_{2}^{p,q}=0,$ we
obtain $p=0.\Box \smallskip $\medskip \newline

\textbf{Remark 3.5. } \textit{Let} $n$ \textit{be an arbitrary positive
integer and} $p,q$ \textit{be two arbitrary integers. Let }$\left(
r_{n}^{p,q}\right) _{n\geq 1}$ \textit{be the generalized Pell-
Fibonacci-Lucas numbers. The following relation is true} 
\begin{equation*}
pP_{n+1}+qQ_{n}=r_{n}^{p,q}+r_{n+1}^{2p,0},\forall ~n\in \mathbb{N}-\{0\}.
\end{equation*}%
\textbf{Proof.} 
\begin{equation*}
pP_{n+1}+qQ_{n}=pP_{n-1}+qQ_{n}+2pP_{n}=r_{n}^{p,q}+2pP_{n}=r_{n}^{p,q}+r_{n+1}^{2p,o}.
\end{equation*}%
$\Box \smallskip $\medskip \newline

\textbf{Definition 3.6. } A subring $O\subseteq \mathbb{H}(\alpha ,\beta ) $
is an \textit{order} in $\mathbb{H}(\alpha ,\beta )$ if $O$ is a finitely
generated $\mathbb{Z}$-submodule of $\mathbb{H}(\alpha ,\beta )$ ( see [Vo;
15]).\medskip

\textbf{Proposition 3.7. } \textit{Let} $O$ \textit{be} \textit{the set} 
\begin{equation*}
O=\left\{ \sum\limits_{i=1}^{n}8R_{n_{i}}^{p_{i},q_{i}}|n\in \mathbb{N}%
^{\ast },p_{i},q_{i}\in \mathbb{Z},(\forall )i=\overline{1,n}\right\} \cup
\left\{ 1\right\} .
\end{equation*}%
\textit{This set} \textit{is an order of the quaternion algebra} $\mathbb{H}%
_{\mathbb{Q}}\left( \alpha ,\beta \right) .$ \newline
\smallskip \newline
\textbf{Proof.} We remark that $0$$\in $$O$ (according to Remark 3.4). We
will show that $O$ is a $\mathbb{Z}-$ submodule of $\mathbb{H}_{\mathbb{Q}%
}\left( \alpha ,\beta \right) .$ Indeed, if $n,m\in \mathbb{N}^{\ast },$ $%
a,b,p,q,p^{^{\prime }},q^{^{\prime }}\in \mathbb{Z},$ it results immediately
that 
\begin{equation*}
ar_{n}^{p,q}+br_{m}^{p^{^{\prime }},q^{^{\prime
}}}=r_{n}^{ap,aq}+r_{m}^{bp^{^{\prime }},bq^{^{\prime }}},
\end{equation*}%
which implies the following relation 
\begin{equation*}
aR_{n}^{p,q}+bR_{m}^{p^{^{\prime }},q^{^{\prime
}}}=R_{n}^{ap,aq}+R_{m}^{bp^{^{\prime }},bq^{^{\prime }}}.
\end{equation*}%
So, we obtain that $O$ is a free $\mathbb{Z}-$ submodule of rank $4$ for the
quaternion algebra $\mathbb{H}_{\mathbb{Q}}\left( \alpha ,\beta \right) .$ 
\newline
To obtain an order, we will prove that $O$ is a subring of $\mathbb{H}_{%
\mathbb{Q}}\left( \alpha ,\beta \right) .$ Let $m,n$ be two integers, $n<m.$
We calculate 
\begin{equation*}
8r_{n}^{p,q}\cdot 8r_{m}^{p^{^{\prime }},q^{^{\prime }}}=8\left(
pP_{n-1}+qQ_{n}\right) \cdot 8\left( p^{^{\prime }}P_{m-1}+q^{^{\prime
}}Q_{m}\right) =
\end{equation*}%
\begin{equation}
=64pp^{^{\prime }}P_{n-1}P_{m-1}+64pq^{^{\prime }}P_{n-1}Q_{m}+64p^{^{\prime
}}qP_{m-1}Q_{n}+64qq^{^{\prime }}Q_{n}Q_{m}.  \tag{3.5}
\end{equation}%
Using Proposition 3.3 (i, ii, iii, iv), Remark 3.5 and the equality (3.5),
it results: 
\begin{equation*}
8r_{n}^{p,q}\cdot 8r_{m}^{p^{^{\prime }},q^{^{\prime }}}=8pp^{^{\prime }}
\left[ Q_{n+m-2}+\left( -1\right) ^{n}Q_{m-n}\right] +64pq^{^{\prime }}\left[
P_{m+n-1}+\left( -1\right) ^{n}P_{m+1-n}\right] +
\end{equation*}%
\begin{equation*}
+64p^{^{\prime }}q\left[ P_{m+n-1}+\left( -1\right) ^{n}P_{m-1-n}\right]
+64qq^{^{\prime }}\left[ Q_{n+m}+\left( -1\right) ^{n}Q_{m-n}\right] =
\end{equation*}%
\begin{equation*}
=64\left[ pq^{^{\prime }}P_{m+n-1}+qq^{^{\prime }}Q_{m+n}\right] +64\left[
\left( -1\right) ^{n}p^{^{\prime }}qP_{m-n-1}+qq^{^{\prime }}Q_{m-n}\right] +
\end{equation*}%
\begin{equation*}
+8\left[ \left( -1\right) ^{n}8pq^{^{\prime }}P_{m-n+1}+\left( -1\right)
^{n}pp^{^{\prime }}Q_{m-n}\right] +8\left[ 8p^{^{\prime
}}qP_{m+n-1}+pp^{^{\prime }}Q_{m+n-2}\right] =
\end{equation*}%
\begin{equation*}
=8r_{m+n}^{8pq^{^{\prime }},8qq^{^{\prime }}}+8r_{m-n}^{8\left( -1\right)
^{n}p^{^{\prime }}q,8qq^{^{\prime }}}+8r_{m-n}^{\left( -1\right)
^{n}8pq^{^{\prime }},\left( -1\right) ^{n}pp^{^{\prime
}}}+8r_{m-n+1}^{\left( -1\right) ^{n}16pq^{^{\prime
}},0}+8r_{m+n-2}^{8p^{^{\prime }}q,pp^{^{\prime
}}}+8r_{m+n-1}^{16p^{^{\prime }}q,0}.
\end{equation*}%
We get that $8r_{n}^{p,q}\cdot 8r_{m}^{p^{^{\prime }},q^{^{\prime }}}$$\in $$%
O.$ From here we obtain that $O$ is an order of the quaternion algebra $%
\mathbb{H}_{\mathbb{Q}}\left( \alpha ,\beta \right) .\Box \smallskip $%
\medskip \newline

\textbf{Conclusions.} In this paper, we provided applications of \ several
cases of difference equation of degree three. One of these applications is
in the Coding Theory and we used the obtained numbers to built, in some
special situations, cyclic codes with good properties. We have defined the
generalized Pell-Fibonacci-Lucas quaternion elements and we proved that,
using these elements, we can get a set which is \textit{\ }an order of the
generalized quaternion algebra $\mathbb{H}_{\mathbb{Q}}\left( \alpha ,\beta
\right) .$

Using these approaches, the above obtained results can constitute the start
for a further research in which we intend to study properties and
applications of other difference equations.\bigskip

\textbf{Acknowledgments.} The authors thank referees for their suggestions
and remarks which helped us to improve this paper.

\begin{equation*}
\end{equation*}%
\textbf{References}\newline
\begin{equation*}
\end{equation*}%
\newline
[Ba, Pr; 09] M. Basu, B. Prasad, \textit{The generalized relations among the
code elements for Fibonacci coding theory}, Chaos, Solitons and Fractals,
41(2009), 2517-2525.\newline
[Fl, Sa; 15] C. Flaut, D. Savin, \textit{Quaternion Algebras and Generalized
Fibonacci-Lucas Quaternions}, Adv. Appl. Clifford Algebras, 25(4)(2015),
853-862.\newline
[Fl, Sh;13] C. Flaut, V. Shpakivskyi, \textit{On Generalized Fibonacci
Quaternions and Fibonacci-Narayana Quaternions}, Adv. Appl. Clifford
Algebras, 23(3)(2013), 673--688.\newline
[Fib.] http://www.maths.surrey.ac.uk/hosted-sites/R.Knott/Fibonacci/fib.html%
\newline
[Ko, Oz, Si; 17] M. E. Koroglu, I. Ozbek, I. Siap, \textit{Optimal Codes
from Fibonacci Polynomials and Secret Sharing Schemes}, Arab. J. Math, 2017,
1-12, DOI 10.1007/s40065-017-0171-7,
https://link.springer.com/content/pdf/10.1007\%2Fs40065-017-0171-7.pdf 
\newline
[Li, Xi; 04] S. Ling, C. Xing, Coding Theory-A first course, Cambridge
University Press, 2004.\newline
[Sa; 17] D. Savin, \textit{About Special Elements in Quaternion Algebras
Over Finite Fields}, Adv. Appl. Clifford Algebras, 27(2)(2017), 1801-1813.%
\newline
[St; 06] A.P. Stakhov, \textit{Fibonacci matrices, a generalization of the
\textquotedblleft Cassini formula\textquotedblright , and a new coding theory%
}, Chaos, Solitons and Fractals, 30(2006), 56-66.\newline
[Vo; 15] J.Voight, \textit{The arithmetic of quaternion algebras, }2015,
available on the author's website:
http://www.math.dartmouth.edu/jvoight/crmquat/book/quatmodforms-041310.

\begin{equation*}
\end{equation*}

\bigskip

Cristina FLAUT

{\small Faculty of Mathematics and Computer Science, Ovidius University,}

{\small Bd. Mamaia 124, 900527, CONSTANTA, ROMANIA}

{\small http://cristinaflaut.wikispaces.com/;
http://www.univ-ovidius.ro/math/}

{\small e-mail: cflaut@univ-ovidius.ro; cristina\_flaut@yahoo.com}

\medskip \medskip \qquad\ \qquad\ \ 

Diana SAVIN

{\small Faculty of Mathematics and Computer Science, }

{\small Ovidius University, }

{\small Bd. Mamaia 124, 900527, CONSTANTA, ROMANIA }

{\small http://www.univ-ovidius.ro/math/}

{\small e-mail: \ savin.diana@univ-ovidius.ro, \ dianet72@yahoo.com}\bigskip
\bigskip

\end{document}